\documentclass[12pt,english]{amsart}
\usepackage{times}
\usepackage[T1]{fontenc}
\usepackage[latin1]{inputenc}
\usepackage{babel}

\makeatletter

\providecommand{\LyX}{L\kern-.1667em\lower.25em\hbox{Y}\kern-.125emX\@}

 \theoremstyle{plain}    
 \newtheorem{thm}{Theorem}[section]
 \numberwithin{equation}{section} 
 \numberwithin{figure}{section} 
 \theoremstyle{plain}    
 \newtheorem{cor}[thm]{Corollary} 
 \theoremstyle{plain}    
 \newtheorem{lem}[thm]{Lemma} 
 \theoremstyle{plain}    
 \newtheorem{prop}[thm]{Proposition} 
 \theoremstyle{definition}
 \newtheorem{defn}[thm]{Definition}
 \theoremstyle{remark}
 \newtheorem{rem}[thm]{Remark}
 \theoremstyle{remark}    
 \newtheorem*{acknowledgement*}{Acknowledgement} 

\usepackage{amssymb}
\usepackage{mathrsfs}
\ifx\undefined\mathscr \else
\let\mathcal\mathscr
\fi

\newcommand{\Inn}{\operatorname{Inn}}
\newcommand{\Out}{\operatorname{Out}}
\newcommand{\Aut}{\operatorname{Aut}}
\newcommand{\Ad}{\operatorname{Ad}}
\newcommand{\Tr}{\operatorname{Tr}}
\newcommand{\id}{\operatorname{id}}

\makeatother
\begin{document}

\title{On the Classification of Full Factors of Type III.}

\author{Dimitri Shlyakhtenko}
\thanks{Research supported by the National
    Science Foundation}

\address{Department of Mathematics, UCLA, Los Angeles, CA 90095, USA}

\email{shlyakht@math.ucla.edu}

\begin{abstract}
We introduce a new invariant \( \mathcal{S}(M) \) for type III factors
\( M \) with no almost-periodic weights. We compute this invariant
for certain free Araki-Woods factors. We show that Connes' invariant
\( \tau  \) cannot distinguish all isomorphism classes of free Araki-Woods
factors. We show that there exists a continuum of mutually non-isomorphic
free Araki-Woods factors, each without almost-periodic weights.
\end{abstract}

\date{\today}

\maketitle

\section{Introduction.}

The necessity to find effective invariants to distinguish full type
III factors comes from the problem of classifying type III factors
naturally occurring in free probability theory of Voiculescu \cite{dvv:book}.
These factors arise as free products of finite-dimensional or hyperfinite
von Neumann algebras \cite{barnett, radulescu1, radulescu2, dykema:typeiii}
and more generally from a second-quantization procedure involving
the free Gaussian functor (the so-called free Araki-Woods factors,
\cite{shlyakht:quasifree:big}). The classification results so far
all relied on Connes' S and Sd invariants \cite{connes, connes:full},
and worked well for factors having almost periodic weights (for example,
in \cite{shlyakht:quasifree:big} a complete classification of free
Araki-Woods factors for which the free quasi-free state is almost-periodic
was given). However, not all free Araki-Woods factors have almost
periodic weights \cite{shlyakht:semicirc}, and the question of their
complete classification remains open.

This paper is devoted to the exploration of free Araki-Woods factors,
having \emph{no} almost-periodic states or weights. The most refined
existing invariant for such factors is Connes' \( \tau  \) invariant
\cite{connes, connes:full}, which measures the {}``degree of continuity''
of the states on \( M \).

We introduce a new invariant for a factor \( M \), given by the intersection
over the set of all normal faithful states (or weights) \( \phi  \)
on \( M \) of the collections of measures, absolutely continuous
with respect to the spectral measure of the modular operator of \( \phi  \).
This invariant is in spirit related to the Connes' Sd invariant, where
the intersection is taken over all \emph{almost-periodic} weights
on \( M \). Amazingly, it turns out that our invariant can be computed
for certain free Araki-Woods classes. 

Using this invariant we are able to produce a pair of non-isomorphic
free Araki-Woods factors, which cannot be distinguished by their \( \tau  \)
invariant.

\begin{acknowledgement*}
I would like
to thank MSRI and the organizers of the operator algebras program
for their hospitality and for the encouraging atmosphere. I am also
grateful to Professors N. Brown, T. Giordano, U. Haagerup, Y. Ueda
and D.-V. Voiculescu for many useful discussions.
\end{acknowledgement*}

\section{Some examples of topologies associated to unitary representations
of \protect\( \mathbb {R}\protect \)}

The purpose of this section is to set notation and to prove certain
results about unitary representation of \( \mathbb {R} \), which
are needed in the rest of the paper.

\subsection{Spectral measures of group representations.}

Let \( \mu  \) be a measure on a measure space \( X \). Denote by
\( \mathcal{C}_{\mu } \) the collection of measures on \( X \)\[
\mathcal {C}_{\mu }=\{\nu :\mu (Y)=0\Rightarrow \nu (Y)=0\textrm{ for all measurable }Y\subset X\}.\]
In other words, \( \mathcal{C}_{\mu } \) is the collection of all
measures, which are absolutely continuous with respect to \( \mu  \).
We shall abuse notation and write \( \mathcal{C}_{T} \) if \( T \)
is a self-adjoint operator. By this we mean the collection of all
measures, absolutely continuous with respect to the spectral measure
of \( T \) on \( \mathbb {R} \).

It should be mentioned that, following an idea of Mackey, \( \mathcal{C}_{\mu } \)
can be viewed as a kind of replacement for the notion of the support
of \( \mu  \). Indeed, any \( \nu \in \mathcal{C}_{\mu } \) can
be written as \( f\cdot \mu  \) for some function \( f \), which
is completely determined except on a set of \( \mu  \)-measure zero.
One can form intersections of two families \( \mathcal{C}_{\mu } \)
and \( \mathcal{C}_{\nu } \); this operation is to remind us of the
notion of intersection of sets. Similarly, inclusion of \( \mathcal{C}_{\mu } \)
and \( \mathcal{C}_{\nu } \) can be thought of as the inclusion of
the support of \( \mu  \) into that of \( \nu  \). We shall say
that a collection of measures \( \mathcal{C} \) is \emph{supported}
on a set \emph{\( Y \)} if for all \( \nu \in \mathcal{C} \), the
complement of \( Y \) has measure zero. Note also that if \( \mu  \)
is \emph{atomic}, then the knowledge of \( \mathcal{C}_{\mu } \)
is exactly equivalent to the knowledge of the set of atoms of \( \mu  \).

We now recall some basic facts about representations and duality of
locally compact abelian groups (see e.g. \cite{hewitt-ross:abstractHarmonicAnalysis}).
Let \( \sigma :G\to U(H) \) be a \( * \)-strongly continuous representation
of \( G \) on a Hilbert space \( H \). Associated to \( \sigma  \)
it is spectral measure class \( \mathcal{C}_{\sigma } \) on the dual
group \( \hat{G} \). The class \( \mathcal{C}_{\sigma } \) can be defined
as the smallest collection of measures, so that (1) if \( \mu \in \mathcal{C}_{\sigma } \)
and \( \mu ' \) is a.c. with respect to \( \mu  \), then \( \mu '\in \mathcal{C}_{\sigma } \)
and (2) for each \( \xi \in H \), the measure obtained as the Fourier
transform of the positive-definite function \( g\mapsto \langle \xi ,\sigma (g)\xi \rangle  \)
belongs to \( \mathcal{C}_{\sigma } \). If \( H \) is separable,
there is a measure \( \mu  \) in \( \mathcal{C}_{\sigma } \), with
the property that it generates \( \mathcal{C}_{\sigma } \) (i.e.,
\( \mathcal{C}_{\sigma } \) is the smallest collection of measures
satisfying (1) and containing \( \mu  \)). We sometimes refer to
this \( \mu  \) as {}``the'' spectral measure of \( \sigma  \).
In particular, \( H \) can be decomposed as \( H=\int _{\chi \in \hat{G}}H_{\chi }d\mu (\chi ) \),
so that \( \sigma =\int _{\chi \in \hat{G}}\sigma _{\chi } \), where
\( \sigma _{\chi }(g)\cdot \xi =\chi (g)\cdot \xi  \), \( \xi \in H_{\mu } \).

If the group \( G \) contains \( \mathbb {R} \) as a dense subgroup,
the representation \( \sigma  \) can be restricted to \( \mathbb {R}\subset G \).
In fact, all representations \( \sigma  \) of \( G \) arise as extensions
of representations of \( \mathbb {R} \), which are continuous not
just in the topology on \( \mathbb {R} \), but also in the restriction
of the topology \( \tau _{G} \) on \( G \) to \( \mathbb {R}\subset G \).
The measure class \( \mathcal{C}_{\sigma } \) and the spectral measure
\( \mu  \) of \( \sigma  \), when interpreted as objects on \( \widehat{\mathbb {R}}\supset \hat{G} \)
become exactly the measure class and the spectral measure of the restriction
of \( \sigma  \) to \( \mathbb {R} \), viewed as a representation
of \( \mathbb {R} \) (this is evident from the direct integral decomposition
formula stated above). In particular, a representation \( \pi  \)
of \( \mathbb {R} \) extends to a representation of \( G \) iff
its spectral measure \( \mathcal{C}_{\pi } \) is supported on \( \hat{G}\subset \widehat{\mathbb {R}} \).

It is customary to choose a particular spectral measure of a representation
of \( \mathbb {R} \) on \( H \), by finding on \( H \) a non-negative
operator \( A \), for which \( \pi (t)=A^{it} \), and letting \( \mu  \)
be the spectral measure of \( A \) (composed with some faithful state
on \( B(H) \)). In particular, denoting by \( \sigma  \) the extension
of \( \pi  \) to \( G \), we have \( \mathcal{C}_{\sigma }=\mathcal{C}_{\pi }=\mathcal{C}_{\Delta } \).

\subsection{Topologies induced by unitary or orthogonal representations of \protect\( \mathbb {R}\protect \).}

Let \( \mu  \) be a measure on the real line, so that \( \mu (-X)=\mu (X) \).
Denote by \( \pi  \) the associated (real or complex) representation
of \( \mathbb {R} \) on \( L^{2}(\mathbb {R},\mu ) \) given by the
map\[
\pi (t)=M_{\exp (2\pi itx)},\]
where \( M_{g} \) denotes the operator of multiplication by \( g \).
Write \( \tau (\mu ) \) for the weakest topology making the map \( t\mapsto \pi (t)\in U(L^{2}(\mathbb {R},\mu )) \)
continuous with respect to the strong operator topology on the unitary
group of \( L^{2}(\mathbb {R},\mu ) \). If \( \mu  \) is not supported
on a cyclic subgroup of \( \mathbb {R} \), \( \pi  \) is injective
and \( \tau (\mu ) \) is a Hausdorff topology.

\begin{prop}
\label{pro:possibleLCcompletions}The completion of \( \mathbb {R} \)
with respect to the topology \( \tau (\mu ) \) is a locally compact
group iff either \( \tau (\mu ) \) is the usual topology on \( \mathbb {R} \),
or \( \mu  \) is atomic (in which case the completion is compact).
\end{prop}
\begin{proof}
Denote by \( (G,\tau ) \) the completion of \( (\mathbb {R},\tau (\mu )) \).
Then \( \mathbb {R}\subset G \) is an inclusion of locally compact
abelian groups; by Pontrjagin duality, this inclusion is dual to the
injective dense inclusion \( \hat{G}\subset \widehat{\mathbb {R}} \). 

By the structure theory for locally compact abelian groups \cite{hewitt-ross:abstractHarmonicAnalysis},
the connected component of identity of \( \hat{G} \) must have the
form \( R\times H \), with \( R\cong \mathbb {R}^{n} \) and \( H \)
compact and connected. 

First note that \( H=\{e\} \). Indeed, the image of \( H \) in \( \widehat{\mathbb {R}} \)
must be a connected compact subgroup of \( \widehat{\mathbb {R}} \),
hence must be the trivial group.

Since there are no continuous injective maps from \( \mathbb {R}^{n} \)
into \( \mathbb {R} \) for \( n>1 \), either \( n=1 \) or \( n=0 \).
If \( n=1 \), all continuous injective group homomorphisms from \( \mathbb {R} \)
to itself are surjective (their image must be a path-connected subgroup
of \( \mathbb {R} \)). Hence injectivity of \( \hat{G}\subset \widehat{\mathbb {R}} \)
requires that \( \hat{G}=R \) in this case, the inclusion be a homeomorphism
onto \( \widehat{\mathbb {R}} \) and hence \( \tau (\mu ) \) be
the usual topology on \( \mathbb {R} \).

If \( n=0 \), \( \hat{G} \) must be discrete. This corresponds to
the completion \( G \) being compact. Moreover, it is not hard to
see that \( \mu  \) must be supported on \( \hat{G}\subset \widehat{\mathbb {R}} \),
since the representation \( \pi  \) with spectral measure \( \mu  \)
must extend (by the definition of \( \tau (\mu ) \)) to the completion
\( G \). Hence \( \mu  \) is atomic.
\end{proof}
\begin{lem}
A sequence \( \{t_{n}\}_{n=1}^{\infty } \) converges to zero in \( \tau (\mu ) \)
iff \( \hat{\mu }(t_{n})\to 1 \), where \( \hat{\mu } \) is the
Fourier transform of \( \mu  \).
\end{lem}
\begin{proof}
Let \( \pi  \) be a representation of \( \mathbb {R} \) associated
to \( \mu  \) as above, and let \( \xi \in L^{2}(\mathbb {R},\mu ) \)
be the constant function \( 1 \). Then\[
\hat{\mu }(t)=\langle \pi (t)\xi ,\xi \rangle .\]
By definition, \( t_{n}\to 0 \) in \( \tau (\mu ) \) iff \( \pi (t_{n})\to 1 \)
strongly. The vector state \( \phi (T)=\langle T\xi ,\xi \rangle  \)
defines a faithful normal state on the commutative von Neumann algebra
\( \pi (\mathbb {R})''\subset B(L^{2}(\mathbb {R},\mu )) \). Hence
strong convergence of \( \pi (t_{n}) \) to \( 1 \) is equivalent
to\[
\Vert \pi (t_{n})-1\Vert _{2}=\phi ((\pi (t_{n})-1)(\pi (t_{n})-1)^{*}))^{1/2}\to 0.\]
In other words, \( \pi (t_{n})\to 1 \) strongly iff\[
\phi (\pi (t_{n})+\pi (t_{n})^{*})\to 1,\]
i.e., \( \Re \hat{\mu }(t_{n})\to 1 \). Since \( |\hat{\mu }(t_{n})|\leq 1 \),
this happens iff \( \hat{\mu }(t_{n})\to 1 \).
\end{proof}
\begin{thm}
\label{thm:continuumTauInv}There exists a continuum of non-atomic
measures \( \mu _{\lambda } \), \( \lambda \in I \) so that the
topologies \( \tau (\mu _{\lambda }) \) are mutually non-equivalent.
\end{thm}
\begin{proof}
Let \( \{c_{n}:n=1,2,\ldots \} \) be a sequence of real numbers,
so that \( c_{n}\geq 0 \), and \( \sum c^{2}_{k}<+\infty  \). Denote
by \( \mu _{n} \) the \( n \)-fold convolution of delta-measures\[
\mu _{n}=(\frac{1}{2}\delta _{c_{1}}+\frac{1}{2}\delta _{-c_{1}})*\cdots *(\frac{1}{2}\delta _{c_{n}}+\frac{1}{2}\delta _{-c_{n}}).\]
Each \( \mu _{n} \) is a symmetric probability measure; and the Fourier
transform of \( \mu _{n} \) is given by\[
\hat{\mu }_{n}(t)=\prod _{k=1}^{n}\cos (2\pi c_{k}t).\]
The Fourier transform of the weak limit \( \mu  \) of \( \mu _{n} \)
is given by the pointwise limit of this expression; the measure \( \mu  \)
is non-atomic (see e.g. \cite{graham-mcgehee}). 

Let \( c_{k}=3^{-k!} \). Let \( 0<\lambda \leq 1 \) be fixed. Let
\( t_{n}=\lambda 3^{n!}=\lambda c_{n}^{-1} \). We claim that \( t_{n}\to 0 \)
in \( \tau (\mu ) \) iff \( \lambda =1 \). Recall that \( t_{n}\to 0 \)
in \( \tau (\mu ) \) iff \( \hat{\mu }(t_{n})\to 1 \).

If \( \lambda <1 \), \( |\hat{\mu }(t_{n})|\leq |\cos (2\pi c_{n}t_{n})|=|\cos (2\pi \lambda )|<1 \),
so that \( t_{n} \) does not converge to \( 0 \) in \( \tau (\mu ) \).

If \( \lambda =1 \),\[
\hat{\mu }(t_{n})=\prod _{k=1}^{n}\cos (2\pi 3^{n!-k!})\cdot \prod _{k>n}\cos (2\pi 3^{n!-k!})=\prod _{k>n}\cos (2\pi 3^{n!-k!}).\]
There exists \( \omega >0 \) so that for \( 0\leq x\leq \omega  \),
\( \cos (2\pi x)\geq 1-49x^{2} \). Hence\[
\hat{\mu }(t_{n})\geq \prod _{k>n}(1-49\cdot 6^{n!-k!})\]
as long as \( k>n \) and \( n \) is such that \( 3^{n!-k!}\leq 3^{n!-(n+1)n!}=3^{-n\cdot n!}<\omega  \). 

Since for \( a\in [0,1/2] \),\[
\lim _{p\to \infty }(1-49a^{p})^{p}=1\]
uniformly, and the function \( p\mapsto (1-49a^{p})^{p} \) is increasing
for any \( a<1 \), given \( 1>\delta >0 \), there exists a \( p=p(\delta )<+\infty  \),
so that\[
(1-49a^{p})>(1-\delta )^{1/p}\]
for all \( a\in [0,1/2] \). Hence letting \( a=\frac{1}{6}, \) we
get that for any \( n \) so that \( k!-n!\geq n\cdot n!>p \),\[
\hat{\mu }(t_{n})\geq \prod _{k>n}(1-\delta )^{1/(k!-n!)}.\]
Hence\[
\log \hat{\mu }(t_{n})\geq \sum _{k>n}\log (1-\delta )\frac{1}{k!-n!}.\]
Since \( \log (1-\delta )<0 \), and\[
\sum _{k<n}\frac{1}{k!-n!}\leq \sum _{k<n}\frac{1}{k!-\frac{k!}{n}}=\sum _{k<n}\frac{1}{k!}\cdot \frac{1}{1-\frac{1}{n}}\leq \sum _{k<n}\frac{1}{k!}<e,\]
we get that\[
\log \hat{\mu }(t_{n})\geq e\log (1-\delta ).\]
Hence\[
\hat{\mu }(t_{n})\geq (1-\delta )^{e}\]
for any \( n \) so that (1) \( 3^{-n\cdot n!}<\omega  \) and (2)
\( n\cdot n!>p(\delta ) \). It follows that \( \hat{\mu }(t_{n})\to 1 \)
as \( n\to \infty  \).

Now fix \( 0<\lambda <1 \) and set\[
\mu _{\lambda }(X)=\mu (\lambda \cdot X)\]
for any Borel set \( X\subset \mathbb {R} \). Then\[
\hat{\mu }_{\lambda }(t)=\hat{\mu }(t/\lambda ).\]
It follows that for any \( 0<\nu \leq \lambda  \), the sequence \( \nu 3^{n!} \)
is convergent to \( 0 \)  in \( \tau (\mu _{\lambda }) \) iff  \( \nu =\lambda  \).
It follows that \( \{\tau (\mu _{\lambda }):0<\lambda \leq 1\} \)
are pairwise non-equivalent.
\end{proof}
The author is indebted to U. Haagerup for communicating to him the
following example:

\begin{thm}
\label{thm:singularUsualTau}There exists a measure \( \mu  \), so
that \( \mu  \) as well as its arbitrary convolution powers \( \mu *\cdots *\mu  \) (any number of
times) are singular with respect to Lebesgue measure, but \( \tau (\mu ) \)
is the usual topology on the additive group of real numbers.
\end{thm}
\begin{proof}
Let \( \mu  \) be as in the proof of Theorem \ref{thm:continuumTauInv},
with \( c_{n}=3^{-n} \). Then\[
\hat{\mu }(t)=\prod _{k\geq 1}\cos (2\pi \frac{t}{3^{k}}).\]
We first claim that \( \tau (\mu ) \) is the usual topology on the
real line. Assume that \( t_{n}\to \infty  \), but \( \hat{\mu }(t_{n})\to 1 \).
Choose \( N \) so that for all \( n>N \), \( t_{n}>9 \). Choose
\( k \) so that \( c=t_{n}3^{-k}>1 \) but \( t_{n}3^{-k-1}=c/3\leq 1 \).
Then\[
\hat{\mu }(t)\leq \cos (2\pi t_{n}3^{-k})\cdot \cos (2\pi t_{n}3^{-k-1})\cdot \cos (2\pi t_{n}3^{-k-2)}=\cos (2\pi c)\cdot \cos (2\pi c/3)\cdot \cos (2\pi c/9),\]
where, \( 1<c\leq 3 \). Let\[
f(c)=\cos (2\pi c)\cdot \cos (2\pi c/3)\cdot \cos (2\pi c/9).\]
It is not hard to see that \( f(c) \) is strictly less than \( 1 \)
on the interval \( 1<c\leq 3 \). It follows that \( \hat{\mu }(t)<1 \)
whenever \( t>9 \). Contradiction.

All convolution powers of \( \mu  \) are singular with respect to
Lebesgue measure (see the discussion of Taylor-Johnson measures \cite{graham-mcgehee}
for examples of similar measures \( \mu  \) but satisfying even stronger
properties than what we need here).
\end{proof}

\subsection{Bimodule decompositions of crossed products.\label{sect:bimodules}}

Let \( (M,\phi ) \) be a von Neumann algebra, \( \phi  \) a faithful
normal state, and let \( G \) be a locally compact abelian group.
Assume that \( \alpha  \) is an action of \( G \) on \( M \), which
leaves \( \phi  \) invariant. Then the crossed product von Neumann
algebra \( C=M\rtimes _{\alpha }G \) contains a canonical copy \( A \)
of the group algebra \( L(G)\cong L^{\infty }(\hat{G}) \); moreover,
the weight \( \phi  \) gives rise to a normal faithful conditional
expectation \( E:C\to A \). Let \( \hat{\psi } \) be a normal faithful
weight on \( A \), and let \( \psi =\hat{\psi }\circ E \). This
is a normal faithful weight on \( C \). Moreover, \( L^{2}(C,\psi ) \)
is an \( A \),\( A \)-bimodule in a natural way.

Fix an isomorphism \( (A,\hat{\psi })\cong L^{\infty }(\hat{G},\nu _{G}) \).
Denote by \( \ell ^{2} \) the Hilbert space with basis \( e_{1},e_{2},\ldots  \)
and by \( \ell _{n}\subset \ell ^{2} \) the subspaces spanned by
\( e_{1},\ldots ,e_{n} \). Given a measure \( \eta  \) on \( X\times X \)
whose projections onto the the coordinate directions on \( X\times X \)
are both equivalent to \( \nu _{G} \), and a multiplicity function
\( n:X\times X\to \mathbb {N}\cup \{\infty \} \), let\[
H(\eta ,n)=L^{2}(X\times X,\eta ,n)\]
be the space of square-integrable functions from \( X\times X\to \ell ^{2} \),
so that \( f(x,y)\in \ell _{n(x,y)} \) for all \( x,y\in X \) (where
for convenience we set \( \ell_{\infty }=\ell ^{2} \)). Endow
\( H(\eta ,n) \) with an \( A \),\( A \)-bimodule structure by
letting\[
(f\cdot h\cdot g)(x,y)=f(x)h(x,y)g(y),\qquad f,g\in A,\quad h\in H.\]
In fact \cite{connes:correspondences, popa:correspondences, connes:ncgeom}
any bimodule over \( A \) can be represented in this way. It is easily
seen that if \( \eta ' \) is another measure on \( X\times X \),
projecting onto \( \mu  \), and \( \eta ' \) is mutually absolutely
continuous with \( \eta  \), then \( H(\eta ,n)\cong H(\eta ',n) \)
as bimodules over \( A \).

Choose vectors \( \xi _{1},\xi _{2},\ldots \in L^{2}(M,\phi ) \),
\( \Vert \xi _{1}\Vert _{2}=1 \), \( \Vert \xi _{i}\Vert _{2}\leq 1 \),
a measure \( \mu  \) on \( \hat{G} \) and a multiplicity function
\( n:\hat{G}\to \mathbb {N} \), so that:\begin{eqnarray*}
\alpha _{g}(\xi _{i})\perp \alpha _{h}(\xi _{j}) &  & \forall i,j\quad \forall g,h\in G\\
L^{2}(M,\phi ) & = & \overline{\textrm{span}}\{\alpha _{g}(\xi _{i}):g\in G,i=1,2,\ldots \}\\
\langle \xi _{i},\alpha _{g}(\xi _{j})\rangle  & = & \hat{\mu }_{i}(g),
\end{eqnarray*}
where \( \mu _{i}=\mu |_{n^{-1}(i)} \) and \( \hat{\cdot } \) denotes
the Fourier transform. Let \( (X,\mu )=(\hat{G},\nu _{G}) \), where
\( \nu _{G} \) is the Haar measure on \( G \). Let\[
\eta (x,y)=\mu (x-y),\qquad n(x,y)=n(x-y)\]
be a measure and a multiplicity function on \( \hat{G}\times \hat{G} \),
and let \( H=H(\eta ,n) \) (note that the projections of \( \eta  \)
onto the coordinate directions are precisely \( \mu *\nu _{G}=\mu (\hat{G})\cdot \nu _{G} \)).
We claim that \( H(\eta ,n)\cong L^{2}(C,\psi ) \) as bimodules.
To see this, one can verify that the map \( p\xi _{i}p\mapsto p(x)p(y)\chi _{n^{-1}(\{i\})} \)
for a projection \( p\in L^{\infty }(\hat{G}) \) with \( \hat{\psi }(p)<+\infty  \),
is a bimodule isometry from the linear span of \( p\xi _{i}p\in L^{2}(C) \)
to \( H(\eta ,\mu ) \).

Note that the measure \( \mu  \) (which is the {}``spectral measure''
of the representation of \( G \) on \( L^{2}(M) \)) is uniquely
determined up to absolute continuity by the \( L(G), L(G) \)-bimodule
structure of \( L^2(C,\psi) \).

\section{Full type III factors.}

Assume that \( M \) is a full factor, so that its group of inner
automorphisms \( \Inn (M) \) is a closed subgroup of the group of
all automorphisms \( \Aut (M) \), endowed with the \( u \)-topology
\cite{haagerup:standard, connes:full}. Let \( \Out (M)=\Aut (M)/\Inn (M) \),
with the quotient topology. Denote by \( \pi  \) the quotient map
from \( \Aut (M) \) to \( \Out (M) \), and by \( \delta  \) the
composition \( \pi \circ \sigma ^{\phi }_{t} \) (which is independent
of \( t \) by Connes' Radon-Nikodym type theorem \cite{connes}).

Assume that the action of \( \mathbb {R} \) on \( M \) by \( t\mapsto \sigma _{t}^{\phi } \)
extends, for some \( \phi  \), to an action of a locally compact
completion \( G \) of \( \mathbb {R} \) (by Proposition \ref{pro:possibleLCcompletions}
above, this means that either \( G \) is just \( \mathbb {R} \),
or \( G \) is compact, and \( \phi  \) is almost-periodic). In this
case, call \( \phi  \) a \( G \)-state (or weight) on \( M \).
Call the crossed product\[
M\rtimes _{\sigma ^{\phi }}G\]
the \( G \)-core of \( M \). It is known \cite{connes, connes:full}
that the \( G \)-core of \( M \) is independent of the choice of the
\( G \)-the state \( \phi  \) (having the property that its modular group
extends to \( G \)).

\begin{defn}
Let \( M \) be full. Then define:\\
(a) \( \tau (M) \) to be the weakest topology on \( \mathbb {R} \)
making the map \( \delta :\mathbb {R}\to \Out (M) \) continuous (this
invariant was introduced by Connes in \cite{connes:full}). Thus \( t_{n}\to 0 \)
in \( \tau (M) \) iff there exists a sequence of faithful normal
semifinite weights \( \phi _{n} \) on \( M \), for which \( \sigma ^{\phi _{n}}_{t_{n}}\to \id  \)
in the \( u \)-topology. \\
(b) Denote by \( \tau _{\phi } \) the weakest topology on \( \mathbb {R} \)
making the modular group \( \sigma _{t}^{\phi } \) of a normal faithful
semifinite weight \( \phi  \) on \( M \) continuous in the \( u \)-topology
on \( \Aut (M) \). Then define \( \overline{\tau }(M) \) to be the
weakest topology among all \( \tau _{\phi } \). In other words, \( t_{n}\to t \)
in \( \overline{\tau }(M) \) iff \( \sigma _{t_{n}-t}^{\phi }\to \id  \)
for some normal faithful semifinite weight \( \phi  \) on \( M \);\\
(c) The \( \mathcal{S} \) invariant to be the intersection\[
\mathcal{S}(M)=\bigcap _{\phi \textrm{ f}.\textrm{n}.\textrm{state on }M}\mathcal{C}_{\bigoplus _{n}(\Delta ^{\phi })^{\otimes n}}.\]

We similarly define\[
\mathcal{W}(M)=\bigcap _{\phi \textrm{ a f}.\textrm{n}.\textrm{ strictly semifinite weight on }M}\mathcal{C}_{\bigoplus _{n}(\Delta ^{\phi })^{\otimes n}}.\]

\end{defn}
Note that \( \mathcal{W}(M)\subset \mathcal{S}(M) \).

Part (c) of the definition is equivalent to the Sd invariant of Connes
\cite{connes:full} if the intersection were to be taken over all
faithful normal \emph{almost-periodic} weights.

Note that since we are dealing only with states in the definition
of \( \mathcal{S}(M) \), the delta measure at \( 1 \) is always
in \( \mathcal{S}(M) \). 

Note that \( \overline{\tau }(M) \) is not weaker than \( \tau (M) \).

\begin{prop}
\label{porp:WGeitherallornothing} \( \mathcal{W}(M) \) is contained
in \( \mathcal{C}_{\lambda } \), the class of the Haar measure on
\( \mathbb {R}_{+}=\widehat{\mathbb {R}} \).
\end{prop}
\begin{proof}
Let \( \phi  \) be any \( G \)-weight on \( M \). Denote by \( \mu  \)
the spectral measure of infinitesimal generator \( \Delta _{\phi } \)
of the unitary group \( \sigma _{t}^{\phi } \) acting on \( (M,\phi ) \).
Denote by \( d \) the function \( x\mapsto x \) on \( \mathbb {R}_{+} \).
On \( B(L^{2}(\hat{R})) \), consider the normal faithful weight \( \psi =\Tr (d\cdot ) \).
Then the modular group of \( \psi  \) is given by conjugation with
\( d^{it} \), and hence \( \psi  \) is a weight on \( B(L^{2}(\widehat{\mathbb {R}})) \).
In particular, the modular group of \( \psi  \), acting on \( L^{2}(\widehat{\mathbb {R}}) \),
is the left regular representation of \( \mathbb {R} \).

Consider the normal faithful weight \( \theta =\phi \otimes \psi  \)
on \( M\otimes B(L^{2}(\widehat{\mathbb {R}}))\cong M \). Then the
spectral measure class \( \Delta _{\phi }\otimes d \) is the Lebesgue
measure \( \lambda  \) on \( \mathbb {R} \), since the left regular
representation of \( \mathbb {R} \) is absorbing for tensor products.
Hence \( \mathcal{W}(M) \) can be no bigger than the class of the
Lebesgue measure.
\end{proof}
Assume that the \( G \)-core of \( M \) is a factor. Note that since
the \( G \)-core has a semifinite normal trace, it is a full iff
its compression by a finite projection is non-\( \Gamma  \). In particular,
if the \( G \) core is a factor and is full, it has no non-trivial
central sequences.

\begin{thm}
Assume that for some \( G \)-state \( \phi  \) on a factor \( M \)
the \( G \)-core is a full factor. Then \( M \) is full.
\end{thm}
\begin{proof}
Let \( C=M\rtimes _{\sigma }G \) be the \( G \)-core of \( M \).
Assume for contradiction that \( M \) is not full. Then by \cite[Corollary 3.6, Proposition 2.8]{connes:full}
there exists a sequence of unitaries \( m_{n}\in M \), \( \phi (m_{n})=0 \),
and so that \( [m_{n},m]\to 0 \) \( * \)-strongly for all \( m\in M \)
and \( [m_{n},\phi ]\to 0 \) in norm on \( M_{*} \) for any \( \phi \in M_{*} \).
View \( m_{n}\in C\supset M \). We claim that \( \{m_{n}\} \) is
a non-trivial central sequence in \( C \). Let \( E:C\to L(G)\subset C \)
be the canonical conditional expectation. Then \( E(m_{n})=\phi (m_{n})=0 \),
so that \( m_{n} \) is not asymptotically scalar. For any \( m\in M \),
\( [m_{n},m]\to 0 \) \( * \)-strongly. Denote by \( U_{g}\in L(G)\subset C \),
\( g\in G \), the implementing unitaries. Then by Connes' results
\cite[Theorem 2.9(3)]{connes:full}, for any \( t\in \mathbb {R}\subset G \),
\( [m_{n},U_{g}]\to 0 \) \( * \)-strongly. Hence \( [m_{n},u]\to 0 \)
\( * \)-strongly for any \( u\in W^{*}(M,U_{t}:t\in \mathbb {R})=C \).
Hence \( m_{n} \) is an asymptotically central sequence in \( C \).
Since \( C \) is assumed to be a full factor, and is of type II\( _{\infty } \),
we have arrived at a contradiction. 
\end{proof}
\begin{thm}
\label{thrm:whenGcoreFull}Assume that for some \( G \)-state \( \phi  \)
on \( M \), the \( G \)-core of \( M \) is a full factor. Then
if for some \( H \) not necessarily locally compact containing \( \mathbb {R} \)
as a dense subgroup, there is an \( H \)-weight \( \psi  \) on \( M \),
one must have that \( H\subset G \). 
\end{thm}
\begin{proof}
Let \( C=M\rtimes _{\sigma ^{\phi }}G \). Let \( L(G)\subset C \)
be the canonical copy of the group algebra of \( G \); for \( g\in G \),
denote by \( w_{g}\in L(G) \) the implementing unitary. We write
\( E_{L(G)} \) for the canonical conditional expectation from \( C \)
onto \( L(G) \).

Assume that \( H\not \subset G \), so that the topology defined by
the inclusion \( \mathbb {R}\subset H \) is not stronger than the
topology defined by the inclusion \( \mathbb {R}\subset G \). Hence
there exists a sequence \( t_{n}\in \mathbb {R} \), so that \( \sigma ^{\psi }_{t_{n}}\to \id  \),
but \( \sigma ^{\phi }_{t_{n}} \) does not converge. Let \( u_{t}=[\phi :\psi ]_{t}\in M\subset C \).
It follows that\[
\Ad _{u_{t}w_{t}}|_{C}=\sigma ^{\psi }_{t}.\]
In particular, \( [u_{t_{n}}w_{t_{n}},x]\to 0 \) \( * \)-strongly
for all \( x\in M\subset C \). Note moreover that \( u_{t}w_{t} \)
commutes with \( u_{s}w_{s} \) (since they form a one-parameter subgroup
of \( U(C) \)). Hence for \( s \) fixed, \( [u_{t_{n}}w_{t_{n}},u_{s}w_{s}]=0 \),
and since \( u_{t_{n}}w_{t_{n}} \) asymptotically commutes with \( M\subset C \),
it follows that also \( [u_{t_{n}}w_{t_{n}},w_{s}]\to 0 \) \( * \)-strongly.
It follows that \( u_{t_{n}}w_{t_{n}} \) is a central sequence in
\( M \). Hence, by the assumption that \( C \) is a full factor,
and by the fact that \( C \) is type II\( _{\infty } \), we find
that \( \lambda _{n}u_{t_{n}}w_{t_{n}}\to 1 \) \( * \)-strongly
for some scalars \( \lambda _{n}\in \mathbb {T} \). But then\[
E_{L(G)}(\lambda _{n}u_{t_{n}}w_{t_{n}})\to 1\]
\( * \)-strongly. Since \( u_{t_{n}}\in M\subset C \), \( E_{L(G)}(u_{t_{n}})\in \mathbb {C} \),
so that \( \lambda '_{n}w_{t_{n}}\to 1 \) \( * \)-strongly for some
sequence \( \lambda '_{n}\in \mathbb{C} \), \( |\lambda|\leq 1 \). Hence \( \lambda'_{n}\pi (t_{n})\to 1 \)
\( * \)-strongly, where \( \pi  \) is the left regular representation
of \( G \) (since the representation of \( G \) on \( L^{2}(C,\Tr ) \)
is a multiple of its left regular representation). 

Choose now \( \phi  \) a function on \( G \), supported in a compact
neighborhood of identity and so that \( \Vert \phi \Vert _{2}=1 \).
Then \( \lambda _{n}\pi (t_{n})\cdot \phi \to \phi  \) in \( L^{2} \).
In particular, it means that the support \( X \) of \( \phi  \)
and its translate \( X+t_{n} \) cannot be disjoint once \( n \)
is sufficiently large. It follows that \( t_{n}\in X-X \) for sufficiently
large \( n \). It follows that \( t_{n}\to 0 \) in \( G \). Contradiction.
\end{proof}
\begin{cor}
\label{thrm:TauBarwhenGcoreFull}If the \( G \)-core of \( M \)
is full, then \( \overline{\tau }(M)=\tau _{G} \), the weakest topology
making the inclusion \( \mathbb {R}\subset G \) continuous.
\end{cor}
We also have:

\begin{prop}
If \( M \) has a \( G \)-state and \( G \subset H \) is a proper
inclusion of locally compact abelian groups,
then the \( H \)-core of \( M \) is not a full factor.
\end{prop}
\begin{proof}
Let \( \phi  \) be a \( G \)-state on \( M \), and denote by \( C \)
the \( H \)-core \( M\rtimes _{\sigma ^{\phi }}H \). Assume that
\( C \) is a factor. By assumption, there exists a sequence \( t_{n}\in \mathbb {R} \),
\( t_{n}\to 0 \) in the topology of \( G \), but \( t_{n} \) not
convergent in \( H \). Denote by \( w_{h}\in C \), \( h\in H \),
the implementing unitaries for the \( H \) action on \( M \). Then
\( \Ad _{w_{t_{n}}}(x)\to x \) \( * \)-strongly for all \( x\in M\subset C \);
moreover, \( \Ad _{w_{t_{n}}}(w)=w \) for all \( w\in W^{*}(w_{h}:h\in H)=L(H) \).
Hence \( w_{t_{n}} \) form a central sequence. Arguing exactly as
in the last paragraph of the proof of Theorem \ref{thrm:whenGcoreFull}
we find that for no sequence of scalars \( \lambda _{n}\in \mathbb {T} \)
does \( \lambda _{n}w_{t_{n}}\to 1 \) \( * \)-strongly in the group
algebra \( L(H) \). Hence \( w_{t_{n}} \) is a non-trivial central
sequence in \( M \).

Choose \( p\in L(H)\subset C \) a projection of finite trace. Then
\( [p,w_{t_{n}}]=0 \) and hence \( pw_{t_{n}}p \) is a central sequence
in \( pCp \), which has a finite trace. Thus \( pCp \) has property
\( \Gamma  \). Hence by Connes' results \cite{connes:full}, \( pCp \)
and hence \( C \) is not full.
\end{proof}

\section{Crossed products, free entropy dimension and the \protect\( \mathcal{S}\protect \)
invariant.}

The main result of this section is a computation of the \( \mathcal{W} \)
and \( \mathcal{S} \) invariants of some type III factors \( M \),
for which the core has a sequence of generators with large free entropy
dimension. We first recall some preliminaries.

\subsection{Free entropy dimension for infinite families of generators.}

It is useful for us to consider Voiculescu's free entropy dimension
in the context of an infinite family of generators \( x_{1},x_{2},\ldots  \)
in a von Neumann algebra \( M \). We point out the necessary modifications
of Voiculescu's approach (\cite{dvv:entropy3,dvv:improvedrandom};
see also \cite{shlyakht:cost:micro}, where a similar modification
was necessary). We freely use the notations of \cite{dvv:entropy2,dvv:entropy3,dvv:improvedrandom}.

Let \( x_{i} \), \( 1\leq i<+\infty  \) and \( y_{i} \), \( 1\leq i<+\infty  \)
be in \( M \). Fix a free ultrafilter \( \omega  \) and an element
\( \kappa  \) in the Stone-Chech compactification of \( (0,1] \),
not lying in this interval. Define\[
\chi ^{\omega }(x_{1},\ldots ,x_{p}:y_{1},y_{2},\ldots ;m,\varepsilon )=\liminf _{q\to \infty }\chi ^{\omega }(x_{1},\ldots ,x_{p}:y_{1},y_{2},\ldots ,y_{q};m,\varepsilon )\]
(note that the \( \liminf  \) is actually a limit in this definition).
Define \( \chi ^{\omega }(x_{1},\ldots ,x_{p}:y_{1},y_{2},\ldots ) \)
in exactly the same way as in \cite{dvv:entropy3}, but using \( \chi  \)
defined above. One still has the property\[
\chi ^{\omega }(x_{1},\ldots ,x_{p}:y_{1},y_{2},\ldots )\leq \chi ^{\omega }(x_{1},\ldots ,x_{p}:z_{1},\ldots ,z_{l})\]
for all \( z_{1},\ldots ,z_{l}\in W^{*}(x_{1},\ldots ,x_{p},y_{1},y_{2},\ldots ) \).

Define\begin{eqnarray*}
 &  & \delta ^{0}_{\omega ,\kappa }(x_{1},x_{2},\ldots ,x_{p}:y_{1},y_{2},\ldots )=\\
 &  & \qquad p-\lim _{\varepsilon \to \kappa }\frac{\chi ^{\omega }(x_{1}+\sqrt{\varepsilon }S_{1},x_{2}+\sqrt{\varepsilon }S_{2},\ldots ,x_{p}+\sqrt{\varepsilon }S_{p}:S_{1},\ldots ,S_{p},y_{1},y_{2},\ldots )}{\log \varepsilon },
\end{eqnarray*}
where \( S_{1},\ldots ,S_{p} \) are a free semicircular family, free
from \( \{x_{i}\}_{i}\cup \{y_{j}\}_{j} \). Now define\[
\underline{\delta }(x_{1},x_{2},\ldots )=\liminf _{p\to \infty }\delta ^{0}_{\omega ,\kappa }(x_{1},\ldots ,x_{p}:x_{p+1},x_{p+2},\ldots ).\]
For a finite family \( x_{1},\ldots ,x_{n} \) this is exactly Voiculescu's
definition of free entropy dimension. In general, \( \underline{\delta }(x_{1},x_{2},\ldots )\leq \liminf _{p\to \infty }\delta ^{0}_{\omega ,\kappa }(x_{1},\ldots ,x_{p}) \).
Moreover,\[
0\leq \underline{\delta }(x_{1},x_{2},\ldots )\]
iff \( W^{*}(x_{1},x_{2},\ldots ) \) can be embedded into \( R^{\omega } \),
the ultrapower of the hyperfinite II\( _{1} \) factor. 

If \( x_{1},\ldots ,x_{p} \) are free form \( x_{p+1},x_{p+2},\ldots  \),
then \( \delta ^{0}_{\omega ,k}(x_{1},\ldots ,x_{p}:x_{p+1},\ldots )=\delta ^{0}_{\omega ,\kappa }(x_{1},\ldots ,x_{p}) \).
In particular, if the families \( \{x_{1}\},\ldots ,\{x_{p}\},\ldots ,\{y_{1},y_{2},\ldots \} \)
are free, and \( \{y_{1},y_{2},\ldots \}'' \) is embeddable, we get
by \cite{dvv:improvedrandom} that\begin{eqnarray*}
\underline{\delta }(x_{1},y_{1},x_{2},y_{2},\ldots ) & = & \lim _{p\to \infty }\delta _{\omega ,k}^{0}(x_{1},\ldots x_{p},y_{1},\ldots ,y_{p}:y_{1},y_{2},\ldots )\\
 & = & \lim _{p\to \infty }\delta _{\omega ,k}^{0}(x_{1},\ldots ,x_{p})+\underline{\delta }(y_{1},y_{2},\ldots )\\
 & \geq  & \sum _{k}\delta (x_{k}).
\end{eqnarray*}

\begin{defn}
Let \( M \) be a II\( _{1} \) von Neumann algebra. Denote by\[
\delta (M)=\sup _{x_{1},x_{2},\ldots \in M}\underline{\delta }(x_{1},x_{2},\ldots )\]
where the supremum is taken over all self-adjoint families (finite
or infinite) \( x_{1},x_{2},\ldots  \) of generators of \( M \).
If \( N \) is type II\( _{\infty } \), we write \( \delta (N) \)
for the supremum over all finite-trace projections \( p\in N \) of
\( \delta (pNp) \).
\end{defn}
\begin{rem}
It is quite likely that \( \delta (M)\in \{-\infty,0,1,+\infty \} \) if \( M \)
is type II\( _{\infty } \). Note also that \( \delta (M)\leq \delta (M\otimes B(H)) \)
for all \( M \). 
\end{rem}
While \( \delta (M) \) is clearly an invariant of \( M \), our inability
to prove that the number \( \delta (x_{1},x_{2},\ldots ) \) is independent
of the choice of generators \( x_{1},x_{2},\ldots  \) \cite{dvv:entropy2, dvv:entropy3, dvv:improvedrandom}
results in the inability to compute the exact value of \( \delta  \)
for infinite-dimensional von Neumann algebras. However, as we pointed
out above, if \( M=L(\mathbb {F}_{n})*N \), with \( N\subset R^{\omega } \)
and \( n=1,2,\ldots  \) or \( +\infty  \), we have \( \delta (M)\geq n \)
and \( \delta (M\otimes B(H))\geq n \) (in fact, \( =+\infty  \)
by \cite{dykema-radulescu:compressions}). Furthermore, as Voiculescu
proved in \cite{dvv:entropy3}, \( \delta (R)=1 \) if \( R \) is
the hyperfinite II\( _{1} \) (or II\( _{\infty } \)) factor; more
generally, \( \delta (M)\leq 1 \) if \( M \) has property \( \Gamma  \)
or has a Cartan subalgebra (since these properties are inherited by
compressions of a von Neumann algebra, these statements are valid
for \( M \) type II\( _{1} \) or type II\( _{\infty } \)).

The following theorem essentially follows from the results of \cite{dvv:entropy3};
we sketch the necessary modifications of the proof coming from the
fact that we may be dealing with infinite families of generators.

\begin{thm}
\label{thrm:cannotbedisjoint}Let \( M \) be a II\( _{1} \) or II\( _{\infty } \)
von Neumann algebra. Let \( L^{\infty }(X,\mu )\cong A\subset M \)
be a diffuse abelian subalgebra, so that \( \Tr _{M}|_{A} \) is semifinite.
View \( L^{2}(M) \) as an \( A,A \) bimodule, and identify it with
\( H(\eta ,n) \) for some measure \( \eta  \) on \( X\times X \).
Assume that \( \eta  \) is disjoint from \( \mu \times \mu  \),
i.e., \( X\times X=Y_{1}\cup Y_{2} \), so that \( \eta (Y_{1})=0 \)
and \( (\mu \times \mu )Y_{2}=0 \). Then \( \underline{\delta }(M)\leq 1 \).
\end{thm}
\begin{proof} (Sketch). 
We first reduce to the case that \( M \) is type II\( _{1} \). Given
\( t\in (0,+\infty ) \), let \( p\in A \) be a finite projection,
corresponding to the characteristic function of some set \( Y\subset X \),
\( \mu (Y)<+\infty  \). Then view \( pMp \) as a bimodule over \( pAp \).
It is not hard to see that \( pMp \) can be identified with \( H(\eta ',n') \),
with \( \eta ' \) absolutely continuous with respect to \( \eta |_{Y\times Y\subset X\times X} \),
\( n'=n|_{Y\times Y\subset X\times X} \). If \( \eta  \) is disjoint
from \( \mu \times \mu  \), then \( \eta ' \) is disjoint from \( \mu '=\mu |_{Y} \).
If the statement of the theorem can be proved for \( pMp \) and \( pAp\subset pMp \),
we would have that \( \delta (qMq)\leq 1 \) for all \( q\in M \)
of finite trace (since \( qMq \) depends up to isomorphism only on
the trace of \( q \)). Hence by definition we get that \( \delta (M)\leq 1 \).

Let \( x_{1},x_{2},\ldots \in M \) be a sequence of generators of
\( M \); by rescaling (which does not affect \( \underline{\delta }(x_{1},x_{2},\ldots ) \)),
assume that \( \Vert x_{j}\Vert =1 \). By the hypothesis, given \( \omega ,\delta >0 \)
and a measure \( \eta ' \) in the absolute continuity class of \( \eta  \),
there exists an \( N=N(\eta ',\omega ,\delta ) \) and a finite family
of \( N \) disjoint measurable subsets \( X_{i},i\in I \) of \( X \),
each of measure \( 1/N \), a subset \( K\subset I\times I \), so
that \( X=\bigcup X_{i} \) and \( \eta '(Y_{2}\setminus \bigcup _{(i,j)\in K}X_{i}\times X_{j})<\omega  \),
\( (\mu \times \mu )(\bigcup _{(i,j)\in K}X_{i}\times X_{j})=|K|/N^{2}<\delta  \).
It follows that for each fixed \( T \), \( \delta  \) and \( \omega  \),
there are projections \( p_{1},\ldots ,p_{N}\in A \) of trace \( 1/N \)
(corresponding to the characteristic functions of \( X_{1},\ldots ,X_{N} \)
in the identification \( A\cong L^{\infty }(X,\mu ) \)), so that:\[
\Vert x_{t}-\sum _{i,j\in K}p_{i}x_{t}p_{j}\Vert _{2}<\omega \]
and\[
\frac{|K|}{N^{2}}<\delta \]
for all \( 1\leq t\leq T \). Using Voiculescu's result \cite{dvv:entropy3}
and the fact that \( p_{1},\dots ,p_{n}\in W^{*}(x_{1}+\sqrt{\varepsilon }S_{1},\ldots ,x_{T}+\sqrt{\varepsilon }S_{T},S_{1},\ldots ,S_{T},x_{T+1},x_{T+2},\ldots ) \),
we get the estimate\begin{eqnarray*}
\chi (x_{1}+\sqrt{\varepsilon }S_{1},\ldots ,x_{T}+\sqrt{\varepsilon }S_{T}:S_{1},\ldots ,S_{T},x_{T+1},x_{T+2},\ldots ) &  & \\
\leq \chi (x_{1}+\sqrt{\varepsilon }S_{1},\ldots ,x_{T}+\sqrt{\varepsilon }S_{T}:p_{1},\dots ,p_{N}) &  & \\
\leq (T(1-\delta )-1)\log (\varepsilon +\omega )+C &  & 
\end{eqnarray*}
where \( C \) is a constant, independent of \( \omega  \), \( \varepsilon  \)
and \( \delta  \). Letting \( \omega \to 0 \) first, we conclude
that\begin{eqnarray*}
T-\lim _{\varepsilon \to \kappa }\frac{\chi (x_{1}+\sqrt{\varepsilon }S_{1},\ldots ,x_{T}+\sqrt{\varepsilon }S_{T}:S_{1},\ldots ,S_{T},x_{T+1},x_{T+2},\ldots )}{\log \varepsilon } &  & \\
\leq T-T(1-\delta )+1=1+\delta T. &  & 
\end{eqnarray*}
Since \( \delta >0 \) is arbitrary, it follows that \( \delta ^{0}_{\omega ,\delta }(x_{1},\ldots ,x_{T}:x_{T+1},x_{T+2},\ldots )\leq 1 \)
for all \( T \). Hence \( \underline{\delta }(x_{1},x_{2},\ldots )\leq 1 \).
Since the sequence of generators \( \{x_{j}\} \) was arbitrary, we
get that \( \delta (M)\leq 1 \).
\end{proof}
In a similar way one sees that \( \delta (M)>1 \) implies that \( M \)
is a non-\( \Gamma  \) factor.

\subsection{Free entropy dimension and crossed products.}

Using Voiculescu's estimates from \cite{dvv:entropy3}, stated above in Theorem
\ref{thrm:cannotbedisjoint}, we record the following theorem (due
to Voiculescu, but formulated by him under the additional hypothesis
that \( M \) be finitely generated):

\begin{thm}
\label{thrm:crossedproductspectrum}Let \( (M,\phi ) \) be a von
Neumann algebra. Let \( \alpha  \) be an action of a locally compact
nondiscrete abelian group \( G \) on \( M \). Assume that \( \alpha  \)
preserves the state \( \phi  \) on \( M \). Denote by \( U_{g}:L^{2}(M,\phi )\to L^{2}(M,\phi ) \)
the unitaries extending \( \alpha (g):M\to M \). Let \( \mu \in M(\hat{G}) \)
be the spectral measure of the representation \( g\mapsto \bigoplus _{n}(U\oplus \bar{U})^{\otimes n}_{g} \)
(here \( \bar{U} \) denotes the conjugate representation). Let \( C=(M\rtimes _{\alpha }G)\otimes B(H) \).
Assume that for some normal faithful weight \( \psi  \) on \( L(G) \),
the composition \( \psi \circ E_{L(G)}:C\to \mathbb {R} \) is a normal
faithful semifinite trace on \( C \). 

Then if \( C \) satisfies \( \delta (C)>1 \), \( \mathcal{C}_{\mu } \)
must contain the Haar measure of \( G \).
\end{thm}
\begin{proof}
View \( C \) as a bimodule over the abelian subalgebra \( L(G,\gamma )\cong L^{\infty }(G,\mu ) \)
of \( M\rtimes _{\alpha }G \). This bimodule can be identified with
\( H(\eta ,n) \) for some measure \( \eta  \) on \( \hat{G}\times \hat{G} \)
and some multiplicity function \( n \) (see \S\ref{sect:bimodules}).

By Theorem \ref{thrm:cannotbedisjoint}, \( \eta  \) cannot be disjoint
from the product measure \( \nu _{\hat{G}}\times \nu _{\hat{G}} \)
on \( \hat{G}\times \hat{G} \). It follows that \( \mu  \) cannot
be disjoint from the Haar measure \( \nu _{\hat{G}} \). We may assume
that \( \mu  \) is symmetric. We claim that \( \nu _{\hat{G}} \)
is absolutely continuous with respect to \( \mu '=\sum _{n\geq 1}\frac{1}{2^{n}}\mu ^{*n} \),
which is the spectral measure of \( g\mapsto \bigoplus _{n}U_{g}^{\otimes n} \).
We must show that if \( \mu '(X)=0 \), then also \( \nu _{\hat{G}}(X)=0 \)
for all Borel subsets \( X\subset \hat{G} \). Assume to the contrary
that \( \nu _{\hat{G}}(X)>0 \) but \( \mu '(X)=0 \). Then \( \mu ^{*n}(X)=0 \)
for all \( n \). Since \( \mu  \) is not disjoint from \( \nu _{G} \),
there exists a subset \( Y \) of \( \hat{G} \), for which \( \nu _{\hat{G}}|_{Y}=f\cdot \mu |_{Y} \).
By modifying \( \mu  \) and \( Y \), we may assume that \( f=1 \)
and \( \mu (Y)=\nu _{\hat{G}}(Y)=1 \). Then \( \nu _{\hat{G}}=f_{n}\cdot \mu ^{*n} \)
on \( nY=Y+Y+\ldots +Y \) (\( n \) times). Moreover, since \( \mu  \)
and \( \nu _{\hat{G}} \) are symmetric, we may assume that \( Y=-Y \).
Hence we find that \( \nu _{\hat{G}}=f\cdot \mu ' \) on the subgroup
of \( \hat{G} \) generated by \( Y \).

Clearly, \( X \) (modulo a set of \( \nu _{G} \)-measure zero) is
contained in the complement of \( H \), and \( \nu _{G}(H)\neq 0 \).
Since \( \nu _{G} \) is \( \sigma  \)-finite, there exists a countable
sequence \( x_{n}\in \hat{G} \), so that (after possibly making \( X \)
smaller by a set of \( \nu _{G} \)-measure zero), \( X\subset \bigcup _{n}(x_{n}+H) \).

Let \( x\in \hat{G} \). Then for any measure \( \sigma  \) on \( \hat{G} \)
of finite total mass,\[
(\mu _{a}*\nu _{G}|_{H})(x)=\sigma (x+H).\]
It follows that if there are finite measures \( \sigma _{n} \), absolutely
continuous with respect to \( \mu ' \) and so that \( \sigma _{n}(x_{n}+H)\neq 0 \),
then \( \mu '*(\sum \frac{1}{2^{n}}\sigma _{n}) \) (and hence \( \mu '*\mu ' \)
and hence \( \mu '\sim \mu '*\mu ' \)) give \( X \) a nonzero measure
(which would be a contradiction). Hence \( \mu '(x_{n}+H)=0 \) for
some \( n \). Let \( p\in L^{\infty }(\hat{G}) \) be the projection
corresponding to the characteristic function of \( H \), and let
\( q\in L^{\infty }(\hat{G}) \) be the projection corresponding to
the characteristic function of \( H+x_{n} \). Let \( (s,t)\in H\times (H+x_{n})\subset \hat{G}\times \hat{G} \).
Then \( s-t\in H-H+x_{n}\in H+x_{n} \). Since \( \mu '(x_{n}+H)=0 \),
it follows that the characteristic function of \( H\times (H+x_{n}) \)
is zero in \( L^{2}(\hat{G}\times \hat{G},\eta ,n) \), where \( \eta  \)
and \( n \) are as in \S\ref{sect:bimodules}. But this implies that
\( pCq=0 \), so that \( p \) and \( q \) are not equivalent in
\( C \). Hence \( C \) is not a factor. Hence \( \delta (C)\leq 1 \).
Contradiction.
\end{proof}
A similar statement holds also for actions preserving a normal faithful
semifinite weight on \( M \).

\subsection{Consequences for the \protect\( \mathcal{S}\protect \) invariant.}

\begin{cor}
\label{corr:spectmeasuremustcontain}Assume that for some normal faithful
\( G \) state \( \phi  \) on \( M \), the \( G \)-core \( C \)
of \( M \) satisfies \( \delta (C)>1 \). Then for any other n.f.s.
weight \( \psi  \) on \( M \), the Haar measure on \( \hat{G} \)
is contained in the spectral measure of the action of \( G \) on
\( \bigoplus _{n}L^{2}(M,\psi )^{\otimes n} \).
\end{cor}
\begin{proof}
This is immediate from \( C=M\rtimes _{\sigma ^{\phi }}G \) and Theorem
\ref{thrm:crossedproductspectrum}.
\end{proof}
\begin{thm}
\label{thrm:deltaandSG}Assume that the core \( C \) of \( M \)
satisfies \( \delta (C)>1 \). Then \( \mathcal{W}(M)=\mathcal{C}_{\lambda } \),
where \( \lambda  \) is Lebesgue measure on the multiplicative group
\( \mathbb {R}_{+} \).\\
If in addition there exists a state \( \phi  \) on \( M \), for
which the spectral measure of the modular group is \( \lambda +\delta _{1} \),
then \( \mathcal{S}(M)=\mathcal{C}_{\lambda +\delta _{1}} \).
\end{thm}
\begin{proof}
By Corollary \ref{corr:spectmeasuremustcontain}, we get that \( \mathcal{C}_{\lambda }\subset \mathcal{W}(M)\subset \mathcal{S}(M) \).
Moreover, \( \mathcal{W}(M)\subset \mathcal{C}_{\lambda } \) in general. 

The proof of the second part of the statement proceeds exactly as
the proof above, with the exception that we now know that \( \psi  \)
can be chosen to be a state, and hence \( \mathcal{S}(M)\subset \mathcal{C}_{\lambda +\delta _{1}} \). 
\end{proof}

\section{Applications to free Araki-Woods factors.}

\subsection{\protect\( G\protect \)-core for certain free Araki-Woods factors.}

Let \( \hat{G}\subset \mathbb {R} \), and denote by \( \sigma  \)
its Haar measure. 

Let \( \nu  \) be a measure on \( \hat{G} \), which is symmetric,
\( \nu (X)=\nu (-X) \). Extending \( \nu  \) to all of \( \mathbb {R} \)
by \( \nu (X)=\nu (X\cap \hat{G}) \) gives us a measure on the real
line. 

Let \( H=L^{2}(\mathbb {R},\nu )=L^{2}(\hat{G},\nu ) \). Denote by
\( H_{\mathbb {R}} \) the subspace of \( H \) consisting of functions
with the property that \( f(x)=\overline{f(-x)} \). Then \( H_{\mathbb {R}} \)
is a real subspace of \( H \), and the restriction of the inner product
on \( H \) to \( H_{\mathbb {R}} \) is real-valued. Moreover, the
one-parameter group of unitary operators\[
U_{t}:t\mapsto \mathcal{M}_{\exp (it)}\]
of multiplication operators acting on \( H \) leaves \( H_{\mathbb {R}} \)
invariant and hence defines an action of \( \mathbb {R} \) on this
real Hilbert space. 

Note that if we consider the dual inclusion \( \mathbb {R}\subset G \),
then the map\[
t\mapsto \mathcal{M}_{\exp (it)}:L^{2}(\hat{G},\nu )\to L^{2}(\hat{G},\nu )\]
extends to the map\[
g\mapsto \mathcal{M}_{\langle g,\cdot \rangle }:L^{2}(\hat{G},\nu )\to L^{2}(\hat{G},\nu )\]
where \( \langle g,\cdot \rangle  \) denotes the function \( \langle g,\cdot \rangle (\chi )=\chi (g) \),
\( g\in G \), \( \chi \in \hat{G} \). Hence \( U_{t} \) extends
to an action \( U_{g} \) of \( G \) on \( H \); it is not hard
to see that once again \( H_{\mathbb {R}} \) is invariant under \( U_{g} \),
\( g\in G \), and hence \( G \) acts on the real Hilbert space \( H_{\mathbb {R}} \)
as well. Note that \( U_{g} \) is isomorphic to the left regular
representation of \( G \). In particular, the spectral measure of
the infinitesimal generator of \( t\mapsto U_{t} \) is \( \nu  \).

Let \( \Gamma (H_{\mathbb {R}},U_{t})'' \) be the free Araki-Woods
factor \cite{shlyakht:quasifree:big} associated to the action \( U_{t} \)
of \( \mathbb {R} \) on \( H_{\mathbb {R}} \), and let \( \phi  \)
denote the free quasi-free state on \( \Gamma (H_{\mathbb {R}},U_{t})'' \).
For convenience we shall write \( \Gamma (\mu ) \) for this von Neumann
algebra.

\begin{thm}
\label{thrm:coreofFAFactors}Let \( \sigma  \) be the Haar measure
on \( G \). Then the \( G \)-core of \( M=\Gamma (\sigma )'' \)
is isomorphic to \( L(\mathbb {F}_{\infty })\otimes B(H) \).
\end{thm}
\begin{proof}
We first note that in the case that \( G \) is compact, the \( G \)-core
is the so-called discrete core of \( M \), and the claimed isomorphism
was already proved by Dykema (\cite{dykema:typeiii}; see \cite{shlyakht:quasifree:big}
for the argument reducing the case of a general Araki-Woods factor
to the form which can utilize Dykema's results). Therefore, we proceed
under the assumption that \( G \) is not compact and hence \( G=\mathbb {R} \),
\( \nu _{\hat{G}} \) is the Lebesgue measure. In particular, \( \nu _{\hat{G}} \)
is non-atomic. 

Let \( C \) denote the core. Let \( \xi \in H_{\mathbb {R}} \) be
a cyclic vector for \( U_{g} \), \( g\in G \) (one can take, for
example, any function \( f(x)=f(-x) \), which is nowhere zero on
\( \hat{G} \), and which lies in \( L^{2}(\hat{G},\sigma ) \)).
Let \( \phi  \) be the positive-definite function on \( G \) associated
to \( \xi  \), \( \phi (g)=\langle \xi ,U_{g}\xi \rangle  \). Let
\( \mu  \) be the Fourier transform of \( \phi  \) (viewed as a
measure on \( \hat{G} \)). By \cite{shlyakht:amalg} and \cite{shlyakht:semicirc},\[
C\cong \Phi (L(G),\eta )\]
where \( \eta  \) is a completely positive map from \( L(G)\cong L^{\infty }(\hat{G})\to L^{\infty }(\hat{G}) \)
given by\[
h\mapsto h*\mu ,\]
\( * \) denoting the convolution on measures on \( \hat{G} \). Notice
that the measure \( \mu  \) is just the measure resulting from applying
the state \( \langle \xi ,\cdot \xi \rangle  \) to the spectral measure
of the infinitesimal generator of \( U_{g} \). Hence \( \mu  \)
is absolutely continuous with respect to the Haar measure \( \sigma  \)
on \( \hat{G} \).

The \( L(G) \), \( L(G) \) bimodule associated to this completely
positive map is\[
L^{2}(\hat{G}\times \hat{G},\hat{\mu }),\]
where \( \hat{\mu }(\chi ,\chi ')=\mu (\chi -\chi ') \), with \( L(G)\cong L^{\infty }(\hat{G}) \)
acting by\[
(f\zeta g)(\chi ,\chi ')=f(\chi )\zeta (\chi ,\chi ')g(\chi ').\]
The real Jordan sub-bimodule of this bimodule (cf. \cite{shlyakht:semicirc})
is generated by the constant function \( 1 \). 

By arguing exactly as in \cite{shlyakht:amalg}, it follows that\[
C\cong \Phi (L^{\infty }(G),\Tr )\cong L(\mathbb {F}_{\infty })\otimes B(H)\]
as claimed.
\end{proof}
The same proof works to show that

\begin{thm}
Let \( M=\Gamma (\sigma ) \) as before and \( \phi  \) be the free
quasi-free state on \( M \). Then the \( G \) core of the \( n \)-fold
free product \( (M,\phi )^{*n} \), for any \( n\geq 1 \) or \( n=+\infty  \),
is \( L(\mathbb {F}_{\infty })\otimes B(H) \). 
\end{thm}
\begin{prop}
Let \( (N,\theta ) \) be a full factor with a \( G \)-state \( \theta  \).
Assume that the \( G \) core of \( N \) can be embedded into \( R^{\omega } \).
Denote by \( C \) the \( G \)-core of \( (N,\theta )*(\Gamma (\nu _{\hat{G}}),\phi ) \).
Then \( \delta (C)=+\infty  \). In particular, \( C \) is full.
\end{prop}
\begin{proof}
Fix \( p\in N\rtimes _{\sigma }G \) a projection of trace \( 1 \).
By \cite{shlyakht:semicirc},\begin{eqnarray*}
C & \cong  & (N\rtimes _{\sigma }G)*_{L(G)}(M\rtimes _{\sigma }G)\\
 & \cong  & (N\rtimes _{\sigma }G)*_{L(G)}\Phi (L(G),\Tr )\cong \Phi (N\rtimes _{\sigma }G,\Tr )\\
 & \cong  & (p(N\rtimes _{\sigma }G)p*L(\mathbb {F}_{\infty }))\otimes B(H),
\end{eqnarray*}
the last isomorphism by \cite{shlyakht-popa:universal} (see also
\cite{dykema-radulescu:compressions}). Since by assumption \( N\rtimes _{\sigma }G \)
is embeddable in \( R^{\omega } \), we get that \( \delta (C)=+\infty  \).
\end{proof}
\begin{thm}
There exists a continuum of mutually non-isomorphic free Araki-Woods
factors, each having no almost-periodic weights.
\end{thm}
\begin{proof}
It was shown in \cite{shlyakht:semicirc} that for each topology \( \tau (\mu ) \)
as discussed above, there exists a free Araki-Woods factor whose \( \tau  \)
invariant is exactly \( \tau (\mu ) \). Moreover, a factor \( M \)
has an almost-periodic weight iff the completion of \( \mathbb {R} \)
with respect to \( \tau (\mu ) \) is compact in that topology (hence
\( \mu  \) is atomic, as then \( \mu  \) is supported on the Pontrjagin
dual \( \Gamma \subset \mathbb {R} \), where \( \mathbb {R}\subset \hat{\Gamma } \)
is the inclusion of \( \mathbb {R} \) into its completion with respect
to \( \tau  \)). By Theorem \ref{thm:continuumTauInv}, there exists
a continuum of mutually non-equivalent topologies \( \tau (\mu _{\lambda }) \),
with \( \mu _{\lambda } \) non-atomic. 
\end{proof}
\begin{thm}
There exist two non-isomorphic free Araki-Woods factors, which cannot
be distinguished by their \( \tau  \) invariant.
\end{thm}
\begin{proof}
Let \( M_{1} \) be the free Araki-Woods factor associated to the
representation of \( \mathbb {R} \) with spectral measure \( \delta _{1}+\delta _{-1}+\lambda  \),
where \( \lambda  \) denotes the Lebesgue measure on the additive
group \( \mathbb {R} \). Then by Theorem \ref{thrm:coreofFAFactors}
and by Corollary \ref{corr:spectmeasuremustcontain},\[
\mathcal{C}(M_{1})\supset \mathcal{C}_{\lambda }.\]
 Moreover, we have that \( \tau (M_{1})=\tau (\delta _{1}+\delta _{-1}+\lambda ) \)
is the usual topology on \( \mathbb {R} \) \cite{shlyakht:semicirc}.

Let \( M_{2} \) be the free Araki-Woods factor associated to the
representation of \( \mathbb {R} \) with spectral measure \( \mu +\delta _{-1}+\delta _{1} \),
where \( \mu  \) is as in Theorem \ref{thm:singularUsualTau}. Then
again \( \tau (M_{2})=\tau (\mu +\delta _{-1}+\delta _{1})=\tau (\mu ) \)
is the usual topology on \( \mathbb {R} \). Thus \( \tau (M_{1})=\tau (M_{2}) \).
However,\[
\mathcal{C}(M_{2})\subset \mathcal{C}_{\sum \frac{1}{2^{n}}(\mu +\delta _{1}+\delta _{-1})^{n}},\]
so that\[
\mathcal{C}(M_{2})\cap \mathcal{C}_{\lambda }=\emptyset .\]
Hence\[
\mathcal{C}(M_{2})\neq \mathcal{C}(M_{1})\]
and so \( M_{1} \) and \( M_{2} \) are not isomorphic. In fact,
going through the proof of Theorem \ref{thrm:crossedproductspectrum},
we see that the cores of \( M_{1} \) and \( M_{2} \) are not isomorphic
(one has \( \delta (M_{1})>1 \), \( \delta (M_{2})\leq 1 \)).
\end{proof}
\bibliographystyle{amsplain}

\providecommand{\bysame}{\leavevmode\hbox to3em{\hrulefill}\thinspace}

\end{document}